\documentclass[reqno]{amsart}
\usepackage{amssymb,amscd,amsxtra}
\usepackage{eucal,pstricks}
\input diagxy

\theoremstyle{plain}

\newtheorem{thm}{Theorem}
\newtheorem{prop}[thm]{Proposition}
\newtheorem{proposition}[thm]{Proposition}
\newtheorem{cor}[thm]{Corollary}
\newtheorem{lem}[thm]{Lemma}
\newtheorem{lemma}[thm]{Lemma}

\theoremstyle{definition}

\newtheorem{rem}[thm]{Remark}

\newenvironment{exam*}{
{\medskip\par\noindent\bf Example. 
}}{\vskip 2ex\par}

\renewenvironment{proof}{
{\medskip\par\noindent\sc Proof.}}{\vskip 2ex\par}

\newcommand{\ga}{\alpha}
\newcommand{\gb}{\beta}
\newcommand{\gc}{\gamma}

\newcommand{\gs}{\sigma}

\newcommand{\frakg}{\mathfrak{g}}
\newcommand{\frakh}{\mathfrak{h}}
\newcommand{\frakm}{\mathfrak{m}}

\newcommand{\frakt}{\mathfrak{t}}

\newcommand{\comp}{\mathrel{\scriptstyle\circ}}

\def\dual#1{#1^{\vee}}

\newcommand{\Ad}{\operatorname{Ad}}

\newcommand{\C}{\mathbb{C}}

\newcommand{\GL}{\operatorname{GL}}

\newcommand{\hatt}{\hat{T}}

\newcommand{\inv}{^{-1}}

\newcommand{\pt}{\operatorname{pt}}
\newcommand{\Q}{\mathbb{Q}}
\newcommand{\R}{\mathbb{R}}

\newcommand{\seq}[2]{#1_{1}, \ldots, #1_{#2}}

\newcommand{\Sym}{\operatorname{Sym}}

\newcommand{\ts}{\tilde{s}}
\newcommand{\tq}{\tilde{q}}
\newcommand{\ty}{\tilde{y}}

\newcommand{\tri}{\bigtriangleup}

\newcommand{\Z}{\mathbb{Z}}

\begin{document}
 \title%[Characteristic Numbers]
{Computing Characteristic Numbers Using Fixed Points 
}
%\medskip
%(Nombres caract\'eristiques d'un espace homog\`ene)}
%  \author{Raoul Bott}
% \address{Department of Mathematics\\
% Harvard University\\
% Cambridge, MA 02138}
% \email{bott@math.harvard.edu}
%\thanks{This work and the first author were supported in part by the 
%National Science Foundation.  
%The author acknowledges the hospitality and support over a 
%period of several years of 
%the Institut Henri Poincar\'e, the \'{E}cole Normale Sup\'erieure, 
%the Institut de Math\'ematiques de Jussieu, and the 
%Minist\`ere de l'Enseignement et de la Recherche, France.}
\author{Loring W. Tu}
\address{Department of Mathematics\\
Tufts University\\
Medford, MA 02155-7049}
\email{loring.tu@tufts.edu}
%\begin{abstract}
%Let $G$ be a compact, connected Lie group with maximal torus $T$, and 
%$H$ a closed subgroup containing $T$.
%We work out the Atiyah-Bott-Berline-Vergne localization 
%formula for the homogeneous space $G/H$ under the natural action of 
%the maximal torus $T$.  
%The computation gives explicit formulas for the ordinary 
%and equivariant characteristic numbers of a homogeneous 
%space.
%%\\
%%
%%\medskip
%%\noindent
%%{\sc R\'esum\'e}.
%%Soit  $G$  un groupe de Lie connexe compact, 
%%muni d'un tore maximal  $T$, et soit  $H$  
%%un sous-groupe qui contient  $T$.  
%% Nous \'etudions la formule de localisation 
%% d'Atiyah-Bott-Berline-Vergne pour l'espace 
%% homogne  $G/H$ sous l'action naturelle du 
%% tore maximal  $T$.  Le calcul donne des 
%% formules explicites pour les nombres 
%% caract\'eristiques ordinaires et 
%% \'equivariantes d'un espace homog\`ene.
% \end{abstract}
\date{August 23, 2011}
\maketitle

\noindent
(This article was published in {\it A Celebration of the Mathematical Legacy of Raoul Bott}, CRM Proceedings and Lecture Notes, vol. 50, American Mathematical Society, Providence, RI, 2010, pp.\ 185--206.
This version differs from the published version
in that Sections 6 and 7 have been switched, as they should be.)

\bigskip

This article has its genesis in a course that Raoul Bott gave at Harvard in the fall of 1996 shortly before his official retirement.
The topic of the course was equivariant cohomology, which is simply the cohomology of a group action.
The course was to culminate in the equivariant localization formula,
discovered by Berline and Vergne, and independently by Atiyah and Bott,
around 1982.
When a manifold has a torus action,
the equivariant localization formula, 
while formulated in equivariant cohomology, 
is a powerful tool for doing calculations in the 
\emph{ordinary} cohomology of the manifold.
It extends and simplifies Bott's work several decades earlier on the relationship between characteristic numbers and the zeros of a vector field (\cite{bott67}). 

In one of the lectures during the course, Raoul Bott considered an
action of a circle on a projective space and computed its equivariant
cohomology from the fixed points using the Borel localization theorem.
After class, he asked me if I could do the same for a homogeneous
space $G/H$ of a compact, connected Lie group $G$ by a closed subgroup
$H$ of maximal rank, under the natural action of a maximal torus.
Unbeknownst to me at the time, and perhaps to him also, this problem
had been solved earlier and is in retrospect not so difficult
(see \cite{arabia} and \cite{brion}, which contain much more than this).
As was his wont, Raoul liked to understand everything \emph{ab initio}
and in his own way.

Using Bott's method, I worked out the equivariant cohomology ring 
of $G/H$ from the fixed points of the torus action.
I saw then that some of the lemmas I developed for this calculation
may be used, in conjunction with the equivariant localization formula,
to calculate the ordinary (as well as the equivariant) characteristic
numbers of $G/H$.

The idea of relating integration of ordinary differential forms to
integration of equivariant differential forms is folklore.
It is implicit in Atiyah--Bott (\cite[Section 8]{atiyah-bott}) and
explicitly stated for the Chow ring in Edidin--Graham
(\cite[Proposition 5, p.~627]{edidin-graham}).
When the fixed points are isolated and the manifold is compact, 
by converting ordinary integration to equivariant integration, 
one can hope to obtain the original integral as a finite sum over the
fixed points of the action.
While the idea is simple, its execution in specific examples is not
necessarily so.
The explicit formulas for the characteristic numbers of $G/H$ obtained
here are involved but beautiful, and serve as an affirmation of the
power and versatility of the equivariant localization formula.
Throughout this project, I met with Raoul many times over a period of
several years.
This article is a testimony to his generosity with his time, ideas,
and friendship, and so I think it is particularly appropriate as a
contribution to a volume on his mathematical legacy.
To make the article more self-contained and in an effort to emulate
Raoul Bott's style, about half of the article is exposition of known
results, for example, the computation
from scratch of the ordinary cohomology rings of $G/T$ and $G/H$.
Although the results are known, 
I have also included the computation of the equivariant cohomology rings of $G/T$
and $G/H$, partly because 
Bott's method of using the Borel localization theorem may be new---at
least I have not seen it in the literature---and may be applicable 
to other situations.
I think of this article as an application of the equivariant
localization formula to one nontrivial example.
It is nonetheless a key example, since every orbit of every action is a homogeneous space.
In the hope that the article
 could be understood 
by a graduate student with a modicum of knowledge of 
equivariant cohomology, perhaps someone who has read 
Bott's short introduction to equivariant cohomology \cite{bott},
I have allowed myself the liberty of being more 
expository than in a typical research article.

Throughout this paper, $H^*(\ )$ stands for singular cohomology
with rational coefficients.
The main technical results are the restriction formulas
(Proposition~\ref{p:restrictgt} and Theorem~\ref{t:restrictgh}) and the Euler
class formulas (Proposition~\ref{p:eulergt} and Theorem~\ref{t:eulergh}),
which allow us to apply the Borel localization theorem and the
equivariant localization formula to compact homogeneous spaces.

Using techniques of symplectic quotients, Shaun Martin 
(\cite[Proposition 7.2]{martin}) has obtained a formula for 
the characteristic numbers of a Grassmannian similar to our 
Proposition~\ref{p:charnogkn}.
Some of the ideas of this paper, in particular that of looking at the
fixed points of the action of $T$ on $G/T$, have parallels in the
Atiyah--Bott proof of the Weyl character formula (\cite{bott88}).

It is a pleasure to thank Aaron W.~Brown for patiently listening to
me and giving me feedback, Jeffrey D.~Carlson
for his careful proofreading and comments, and Jonathan Weitsman for
helpful advice.

\section{Line Bundles on $G/T$ and on $BT$}
\label{s:assobundles}

Let $G$ be a compact, connected Lie group and $T$ a maximal torus  in $G$.
The cohomology classes on the homogeneous space $G/T$ and on the classifying space $BT$ all come from the first Chern class of complex line bundles on these spaces.  Both $G/T$ and $BT$ are base spaces of principal $T$-bundles: $G/T$ is the base space of the principal $T$-bundle $G \rightarrow G/T$ and $BT$ is the base space of a universal principal $T$-bundle $ET \rightarrow BT$.
For this reason, we will first construct complex line bundles on the base space of an arbitrary principal $T$-bundle.

A \emph{character} of the torus $T$ is a multiplicative 
Lie group homomorphism of $T$ into $\C^*$.  Let $\hat{T}$ be the 
group of characters of $T$, with the multiplication of 
characters written additively:  for $\ga, \gb \in \hat{T}$ 
and $t\in T$,
\[
t^{\ga+\gb} := t^{\ga} t^{\gb} = \ga(t) \gb(t).
\]
Denote by $\C_{\gc}$ the complex vector space $\C$ with an 
action of $T$ given by the 
character $\gc\colon T \rightarrow \C^*$.
A character $\ga\colon T \rightarrow \C^*$ has the form
$\ga(t_1, \ldots, t_{\ell}) = t_1^{n_1} \cdots t_{\ell}^{n_{\ell}}$,
where $t_i \in S^1$ and $n_i \in \Z$ (\cite[Proposition 8.1, p.~107]{brocker-tomdieck}). 
So the character group $\hat{T}$ is
isomorphic to $\Z^{\ell}$.

Suppose a torus $T$ of dimension $\ell$ acts freely on the right on a 
topological space $X$ so that $X \rightarrow X/T$ is a principal 
$T$-bundle.  
By the mixing construction, a character $\gc$ on $T$ 
associates to the principal bundle $X \rightarrow X/T$ 
a complex line bundle $L(X/T,\gc)$ over $X/T$:
\[
L_{\gc} := L(X/T,\gc) := X \times_T \C_{\gc} := (X \times \C_{\gc})/T,
\]
where $T$ acts on $X \times \C_{\gc}$ by
\[
(x,v)t= (xt, \gc(t\inv)v).
\]
The equivalence class of $(x,v)$ is denoted $[x,v]$.
It is easy to check that as complex line bundles over $X/T$,
\[
L_{\ga +\gb} \simeq L_{\ga} \otimes L_{\gb}.
\]
Hence, $L_{-\ga} \simeq L_{\ga}^{\spcheck}$, the dual bundle of $L$.

The first Chern class of an associated complex line bundle defines a 
homomorphism of abelian groups
\begin{equation} \label{e:characteristic}
c\colon \hat{T} \rightarrow H^2(X/T), \ \gc \mapsto c_1(L(X/T,\gc)). 
\end{equation}
Let $\Sym(\hat{T})$ be the symmetric algebra with rational
coefficients generated by the additive group $\hat{T}$.  The 
group homomorphism \eqref{e:characteristic} extends to a ring 
homomorphism
\[
c\colon \Sym(\hat{T}) \rightarrow H^*(X/T),
\]
called the \emph{characteristic map} of $X/T$, sometimes 
denoted $c_{X/T}$.

We apply the construction of the associated bundle of a 
character $\gc \in \hat{T}$ to two situations:

(i) The classifying space $BT$.  Applied to the 
universal bundle $ET \rightarrow BT$, this construction yields line 
bundles 
\[
S_{\gc} := L(ET/T, \gc) = L(BT,\gamma)
\]
over $BT$ and cohomology classes $c(\gc) = c_1(S_{\gc})$ 
in $H^2(BT)$.  The characteristic map
\[
c\colon \Sym(\hat{T}) \rightarrow H^*(BT)
\]
is a ring isomorphism, since both $\Sym(\hat{T})$
and $H^*(BT)$ are polynomial rings in 
$\ell$ generators and the generators correspond.  If $\chi_1, \dots, 
\chi_{\ell}$ is a basis for the character group $\hat{T}$ 
and $u_i = c_1(S_{\chi_i})$, then $\Sym (\hat{T})$ is the 
polynomial ring $\Q[\chi_1, \dots, \chi_{\ell}]$, and
\begin{align*}
H^*(BT) &=H^*(BS^1 \times \cdots \times BS^1)\\
&\simeq H^*(BS^1) \otimes_{\Q} \cdots \otimes_{\Q} H^*(BS^1)
 \simeq \Q [ u_1, \dots, u_{\ell}],
\end{align*}
because $BS^1 \simeq \C P^{\infty}$.

(ii) The homogeneous space $G/T$.  Applied to the 
principal $T$-bundle $G \rightarrow G/T$, this construction yields 
line bundles
\[
L_{\gc} := L(G/T, \gc)
\]
over $G/T$.

For each character $\gamma$ on the torus $T$, the relationship of the
line bundle $S_{\gc}$ over $BT$ and the line bundle $L_{\gc}$ over
$G/T$ is as follows.
The universal $G$-bundle $EG\rightarrow BG$ factors through $EG\rightarrow 
EG/T \rightarrow BG$.
The total space $EG$ is a contractible space on which $G$ acts freely.
A fortiori, $EG$ is also a contractible space on which the subgroup
$T$ acts freely.
Hence, $EG = ET$.
It follows that  $(EG)/T= (ET)/T = BT$, so there is a commutative diagram
\[
\bfig
\square(0,300)/^{ (}->`>`>`^{ (}->/<400,300>[G`ET`G/T`BT;```]
\square/`>`>`^{ (}->/<400,300>[G/T`BT`{\rm pt}`BG,;```]
\efig
\]
representing $G/T$ as a fiber of the fiber bundle $BT\rightarrow BG$
and the principal $T$-bundle $G\rightarrow G/T$ as the restriction of the
principal $T$-bundle $ET \rightarrow BT$ from $BT$ to $G/T$.
Hence, the associated bundle $L_{\gc} = G \times_{\gc} \C$ is the
restriction of the associated bundle $S_{\gc} = ET \times_{\gc} \C$
from $BT$ to $G/T$.

\section{The actions of the Weyl group}

The Weyl group of a maximal torus $T$ in the compact,
connected Lie group $G$ is $W = N_G(T)/T$, where $N_G(T)$ 
is the normalizer of $T$ in $G$.  The Weyl group is a finite 
reflection group.

We use $w$ to denote both an element of $W$ and a lift of 
the element to the normalizer $N_G(T)$.
The Weyl group $W$ acts on the character group $\hat{T}$ of 
$T$ by
\[
(w\cdot \gc)(t) = \gc(w\inv t w).
\]
This action induces an action on $\Sym(\hatt)$ as ring 
isomorphisms.
Let $R$ be the polynomial ring
\[
R = \Sym(\hatt) = \Q[\chi_1, \dots, \chi_{\ell}]
\]
and $R^W$ the subring of $W$-invariant polynomials.

If the Lie group $G$ acts on the right on a space $X$, then 
the Weyl group $W$ acts on the right on $X/T$ by
\[
r_w (xT) = (xT) w = xwT.
\]
This action of $W$ on $X/T$ induces by the pullback 
an action of $W$ on the line 
bundles over $X/T$  and also 
on the cohomology ring $H^*(X/T)$:
if $L$ is a line bundle on $X/T$ and $a \in H^*(X/T)$, then
\[
w\cdot L = r_w^* L, \qquad w\cdot a = r_w^* a,
\]
where we use the same notation $r_w^*$ for the pullback of a line
bundle and for the pullback of a cohomology class.

\begin{prop}
The action of the Weyl group $W$ on the associated line 
bundles over $X/T$ is compatible with its action on the 
characters of $T$; more precisely, for $w \in W$ and $\gc 
\in \hatt$, 
\[
r_w^* L(X/T, \gc) \simeq L(X/T, w\cdot \gc).
\] 
\end{prop}

\noindent{\sc Indication of proof}.
Define $\varphi\colon r_w^* L(X/T, \gc) \rightarrow L(X/T, w\cdot \gc)$ by
\[
\varphi(xT, [xw, v])= [x,v].
\]
If we use a different representative $xt$ to represent $xT$, 
then
\[
(xT, [xw, v]) = (xtT, [xtw, v'])
\]
and it is easily verified that $v'=\gc(w\inv t\inv w)v$.
Hence, 
\begin{align*}
\varphi(xtT, [xtw,v']) &= [xt, v'] = [xt, (w\cdot 
\gc)(t\inv)v]\\
&= [x,v].
\end{align*}
This shows that $\varphi$ is well defined.  It has the obvious 
inverse map
\[
\psi ([x,v]) = (xT, [xw, v]). \tag*{\qed}
\]
{\vskip 2ex\par}

\begin{cor} \label{c:charw}
The characteristic map $c\colon \Sym(\hatt) \rightarrow H^*(X/T)$ is 
$W$-equivariant.
\end{cor}

\begin{proof}
Let $\gc$ be a character on the torus $T$.
Writing $L_{\gc}$ instead of $L(X/T,\gc)$, it follows from 
the proposition that
\[
r_w^*c(\gc) = r_w^* c_1(L_{\gc}) = c_1(r_w^*L_{\gc})=
c_1(L_{w\cdot \gc}) = c(w\cdot \gc). \tag*{\qed}
\]
\end{proof}

\begin{cor} \label{c:weylaction}
For $\gc$ a character on $T$ and $w$ an element of the 
Weyl group $W$,
\[
w\cdot c_1 (S_{\gc}) = c_1(S_{w\cdot \gc}).
\]
\end{cor}

\begin{proof}
This is a special case of the preceding corollary with $X=ET$.
\qed\end{proof}

\section{Fiber Bundles with Fiber $G/T$}

Let $G$ be a compact, connected Lie group,
$T$ a maximal torus, and $W= N_G(T)/T$ the Weyl group of $T$ in $G$.
Suppose $G$ acts freely on the right on a topological space $X$ so
that $X \rightarrow X/G$ is a principal $G$-bundle.
Then the natural projection $X/T \rightarrow X/G$ is a fiber bundle
with fiber $G/T$.
We will have frequent occasion to call on the following topological lemma.

\begin{lem}\label{l:gmodtbundle}
The rational cohomology of $X/G$ is the subspace of $W$-invariants of
the rational cohomology of $X/T$:
$H^*(X/G) \simeq H^*(X/T)^W$.
\end{lem}

The proof is based on the following two facts from \cite{hsiang}:

\medskip
\noindent
{\bf Fact 1}.  The compact homogeneous space $G/T$ has a cellular
decomposition into even-dimensional cells indexed by the Weyl group.

This is a consequence of the well-known Bruhat decomposition (see
\cite[p.~35]{hsiang})\footnote{Fact 1 may also be obtained
  from Morse theory; there is a proof in \cite{bott56} for $G$ compact,
  connected, and simply connected.}, using the fact that a compact homogeneous space
$G/T$ has a complex description $G_{\C}/B = G/T$, where $G_{\C}$ is
the complexification of $G$ and $B$ is a Borel subgroup containing
$T$.
It implies that $H^*(G/T)$ vanishes in odd
degrees and that the Euler characteristic of $G/T$ is $\chi(G/T) = |W|$.

\medskip
\noindent
{\bf Fact 2}.  If $N = N_G(T)$ is the normalizer of $T$ in $G$, then $G/N$ is acyclic.

\begin{proof}
The projection $G/T \rightarrow G/N$ is a regular covering map with
group $W=N/T$.  Hence, $H^*(G/N)= H^*(G/T)^W$.  It follows that like
$G/T$,  $G/N$ also has nonzero cohomology classes only in even
degrees.
Since 
\[
\chi(G/N) = \frac{1}{|W|} \chi(G/T) = \frac{1}{|W|} |W| = 1,
\]
$H^0(G/N) = \Q$ and $H^k(G/N) =0$ for $k > 0$.
\qed\end{proof}

\medskip
\noindent
{\sc Proof of Lemma 4}.
Factor $X/T \rightarrow X/G$ into $X/T \rightarrow X/N \rightarrow X/G$.
Because $G/N$ is acyclic, the constant function $1$ defines a global
cohomology class
on $X/N$ that restricts to a generator of cohomology on each fiber $G/N$ of
$X/N\rightarrow X/G$.  By the Leray--Hirsch theorem, $H^*(X/N) \simeq
H^*(X/G)$.

Since $X/T \rightarrow X/N$ is a regular covering with group $W = N/T$,
$H^*(X/N) \simeq H^*(X/T)^W$.
Thus, $H^*(X/G) \simeq H^*(X/T)^W$. \qed

\section{Cohomology Rings of $G/T$ and $G/H$}

Let $EG \rightarrow BG$ and $ET \rightarrow BT$ be universal bundles
for the compact, connected Lie group $G$ and its maximal torus $T$.
Since $ET = EG$, $BG=(EG)/G$ and $BT=(EG)/T$.
So the natural projection $BT \rightarrow BG$ is a fiber bundle with fiber $G/T$.

Choose a basis $\chi_1, \ldots, \chi_{\ell}$ for the character group
$\hat{T}$, let $S_{\chi_i}$ be the associated complex line bundles
over $BT$, and set $u_i = c_1(S_{\chi_i})$.
As noted earlier, if $R = \Sym(\hat{T})$, then
\[
H^*(BT) = \Q [ u_1, \ldots, u_{\ell}] \simeq \Q[ \chi_1, \ldots, \chi_{\ell}
] \simeq R.
\]
By Lemma \ref{l:gmodtbundle}, the cohomology of $BG$ is the subring of
$W$-invariants in $H^*(BT)$:
\[
H^*(BG) = H^*(BT)^W = R^W.
\]

\begin{thm} \label{t:cohomgt}
Let $R_+^W$ be the submodule of $R^W$ generated by all homogeneous elements of positive degree, and
$(R_+^W)$ the ideal in $R$ generated by $R_+^W$.
Then $H^*(G/T) \simeq R/(R_+^W)$.
\end{thm}

\begin{proof}
Consider the spectral sequence of the fiber bundle $BT \rightarrow BG$
with fiber $G/T$.
Since both the base $BG$ and the fiber $BT$ have only even-dimensional
cohomology, all the differentials $d_r$ vanish for $r \ge 2$ and the
spectral sequence
degenerates at the $E_2$ term (\cite[Theorem 14.14 and Theorem 14.18]{bott-tu}).
Therefore,
\[
H^*(BT) = E_{\infty} \simeq E_2 \simeq H^*(BG) \otimes_{\Q} H^*(G/T).
\]
In the picture of $E_2$, $H^*(G/T)$ is the zeroth column and $H^*(BG)$ is the bottom row.
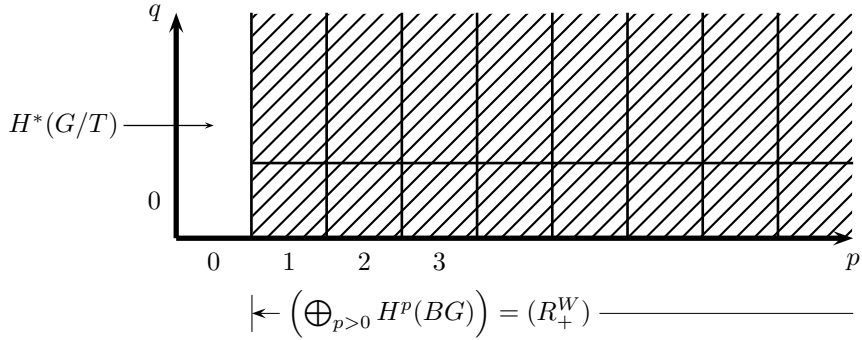
\begin{figure}[h!]
\begin{center}
\begin{pspicture}(-2,-1)(9,3)
\psline[linewidth=2pt]{->}(0,0)(9,0)
\psline[linewidth=2pt]{->}(0,0)(0,3)
\psline[linewidth=1pt](1,0)(1,3)
\psline[linewidth=1pt](2,0)(2,3)
\psline[linewidth=1pt](3,0)(3,3)
\psline[linewidth=1pt](4,0)(4,3)
\psline[linewidth=1pt](5,0)(5,3)
\psline[linewidth=1pt](6,0)(6,3)
\psline[linewidth=1pt](7,0)(7,3)
\psline[linewidth=1pt](8,0)(8,3)
\psline[linewidth=1pt](1,1)(9,1)
\psframe[linestyle=none,fillstyle=hlines](1,0)(9,3)
\rput(-1.5,1.5){$H^*(G/T)$}
\psline[linewidth=.5pt]{->}(-.7,1.5)(.5,1.5)
\uput{.2}[-90](9,0){$p$}
\uput{.2}[180](0,3){$q$}
\uput{.2}[180](0,.5){$0$}
\uput{.2}[-90](.5,0){$0$}
\uput{.2}[-90](1.5,0){$1$}
\uput{.2}[-90](2.5,0){$2$}
\uput{.2}[-90](3.5,0){$3$}
\psset{linewidth=0.2pt,arrowscale=2,tbarsize=7pt}
\psline{|<-}(1,-1)(9,-1)\rput*(3.5,-1){$\left(\bigoplus_{p> 0} H^p(BG) \right)= (R_+^W)$}       
%  \psline{|<-}(2,-2.4)(6,-2.4)\rput*(4,-2.4){$D_2(E_{\infty})$}       
\end{pspicture}
%\caption{The filtration on $E_{\infty}$}
\label{8filtration}
\end{center}
\end{figure}

\noindent
From the picture of $E_2$, one sees that the kernel of the restriction to the fiber:
$H^*(BT) \rightarrow H^*(G/T)$ is 
the shaded area, which is the ideal generated by
$\bigoplus_{p> 0} H^p(BG) = R_+^W$.
Hence,
\begin{equation} \label{e:hgt}
H^*(G/T) \simeq \frac{H^*(BT)}{\left( \bigoplus_{p >0} H^p(BG)\right)}
= \frac{R}{(R_+^W)}.
\end{equation}
This isomorphism is a priori a module isomorphism, but because the
restriction map $H^*(BT) \rightarrow H^*(G/T)$ is a ring homomorphism,
the module isomorphism \eqref{e:hgt} is in fact a ring isomorphism.
\qed\end{proof}

In terms of the basis $\chi_1, \ldots, \chi_{\ell}$ for $\hat{T}$,
let $L_{\chi_i}$ be the line bundle over $G/T$ associated to the
character $\chi_i$ and let $y_i = c_1(L_{\chi_i}) \in H^2(G/T)$.
If $i\colon G/T \rightarrow BT$ is the inclusion map as a fiber of $BT
\rightarrow BG$, then as noted in Section 1,
$L_{\chi_i} = i^* S_{\chi_i}$
and
\[
y_i = c_1(L_{\chi_i}) = c_1(i^*S_{\chi_i}) = i^*c_1(S_{\chi_i}) =
i^*u_i.
\]

Let $H$ be a closed, connected subgroup of maximal rank in the
compact, connected Lie group $G$, and $T \subset H$ a maximal torus.
Denote by $W_G$ and $W_H$ the Weyl groups of $T$ in $G$ and in $H$
respectively.

\begin{thm} \label{t:cohomgh}
If $R=H^*(BT)$, $R^{W_H}$ the subring of $W_H$-invariant elements, and $(R_+^{W_G})$ the ideal generated by $W_G$-invariant elements of positive degree in $R^{W_H}$, then
\[
H^*(G/H) \simeq \frac{R^{W_H}}{\left(R_+^{W_G}\right)}.
\]
\end{thm}

\begin{proof}
The natural projection $G/T \rightarrow G/H$ is a fiber bundle with
fiber $H/T$.
By Lemma \ref{l:gmodtbundle} and Theorem \ref{t:cohomgt},
\[
H^*(G/H) \simeq H^*(G/T)^{W_H}  \simeq \left[ \frac{R}{\left(R_+^{W_G}\right)}\right]^{W_H} = \frac{R^{W_H}}{\left(R_+^{W_G}\right)}. \tag*{\qed}
\]
\end{proof}

\section{Equivariant cohomology and equivariant characteristic classes}

Suppose a topological group $G$ acts on the left on a 
topological space $M$, and 
$EG \rightarrow BG$ is the universal principal $G$-bundle.  Since 
$G$ acts freely on $EG$, the diagonal action of $G$ on $EG 
\times M$,
\[
(e,x)g= (eg, g\inv x),
\]
is also free.  
The space $M_G := EG \times_G M := (EG \times M)/G$ is 
called the \emph{homotopy quotient} of $M$ by $G$, and its 
cohomology $H^*(M_G)$ the \emph{equivariant cohomology} of 
$M$ under the $G$-action, denoted $H_G^*(M)$.
For the basics of equivariant cohomology, we refer to 
\cite{bott} or \cite{guillemin-sternberg}.

A $G$-equivariant vector bundle $E \rightarrow M$ induces a vector 
bundle $E_G \rightarrow M_G$.  An \emph{equivariant characteristic 
class} $c^G(E)$ of $E \rightarrow M$ is defined to be the 
corresponding ordinary 
characteristic class $c(E_G)$ of $E_G \rightarrow M_G$.

By the definition of homotopy quotient, $M_G \rightarrow BG$
is a fiber bundle with fiber $M$.  Let 
$i\colon M \rightarrow M_G$ be the inclusion map as a fiber.
We say that a cohomology class $\tilde{\ga}\in H_G^*(M)$ is 
an \emph{equivariant extension} of $\ga \in H^*(M)$ if 
$i^*\tilde{\ga} = \ga$.  For example, if $E\rightarrow M$ is a 
$G$-equivariant vector bundle, then any characteristic 
class $c(E)$ has an equivariant extension $c^G(E)$, since
\begin{equation} \label{e:extension}
i^*c^G(E) = i^* c(E_G) = c(i^*E_G) = c(E).
\end{equation}

Now suppose $G$ is a compact, connected Lie group with 
maximal torus $T$.  
We associate to a character $\gamma: T \rightarrow \C^{*}$ the complex line 
bundle $L_{\gamma}$ on $G/T$:
\[
L_{\gamma}= G \times_{\gamma}\C.
\]
There is an action of the torus $T$ on $L_{\gamma}$:
\[
t\cdot [g,v] = [tg, v], \quad\text{for } [g,v] \in G\times_{\gamma} \C,
\]
where $[g,v]$ is the equivalence class of $(g,v) \in G \times \C$.
This action is compatible with the action of $T$ on $G/T$ in the sense that 
the diagram
\[
\bfig
\Square[L_{\gamma}`L_{\gamma}`G/T`G/T;t```t]
\efig
\] 
commutes for all $t\in T$.  Therefore, $L_{\gamma} \rightarrow G/T$ is a 
$T$-equivariant line bundle and induces a line bundle $(L_{\gamma})_T$ over 
the homotopy quotient $(G/T)_T$.
Fix a basis $\chi_1, \dots, \chi_{\ell}$ for the 
characters of $T$ and let $\tilde{y}_i = c_1((L_{\chi_i})_T) \in 
H_T^2(G/T)$.  These are equivariant extensions of the 
cohomology classes $y_i := c_1(L_{\chi_i}) \in H^2(G/T)$.

The $T$-equivariant cohomology of a point is
\begin{equation} \label{e:cohombt}
H_T^*(\pt) = H^*(BT) = \Q[u_1, \dots, u_{\ell}], \qquad 
u_i=c_1(S_{\chi_i}).
\end{equation}
For any $T$-space $M$, let $\pi\colon G/T \rightarrow \pt$ be the
constant map.
The pullback map $\pi^*\colon H_T^*(\pt) \rightarrow H_T^*(M)$ makes
$H_T^*(M)$ into an algebra over the polynomial ring $\Q[u_1, \ldots,
  u_{\ell}]$. 

\section{Additive structure of $H_{T}^{*}(G/T)$}

Consider the action of $T$ on $G/T$ by left multiplication.
Write $W = W_G = W_G(T)$ for the Weyl group of $T$ in $G$.
By Theorem \ref{t:cohomgt}, the 
rational cohomology ring of $G/T$ is
\[
H^{*}(G/T) = \Q [\seq{y}{\ell}] / (R_+^W) = R/(R_+^W),
\]
where 
$R=\Q [\seq{y}{\ell}] \simeq H^*(BT)$ and
$(R_+^W)$ is the ideal generated by the $W$-invariant polynomials of 
positive degree in $R$.
In particular, the cohomology $H^{*}(G/T)$ has only even-dimensional cohomology 
classes.  For dimensional reasons, viz., $H^{\text{odd}}(G/T)=
H^{\text{odd}}(BT)=0$, all the differentials $d_{2}, 
d_{3}, \ldots$ in the spectral sequence of the fiber bundle
\begin{equation} \label{e:fibering}
G/T \rightarrow (G/T)_T \rightarrow BT
\end{equation}
vanish, and additively
\begin{align*}
    H^{*}((G/T)_{T}) &= E_{2}\text{-term of the spectral sequence} \\
    &= H^{*}(BT) \otimes_{\Q} H^{*}(G/T) \simeq R\otimes_{\Q} H^*(G/T)
\simeq R \otimes_{\Q}(R/(R_+^W))\\
    &= \Q [\seq{u}{\ell}] \otimes_{\Q} (\Q [\seq{y}{\ell}]/(R_+^W)).
\end{align*}
This shows that $H_{T}^{*}(G/T)$ is a free $R$-module of rank 
equal to $\dim H^{*}(G/T)$.

Moreover, because the differentials $d_r, r\ge 2$, in the spectral
sequence all  
vanish, the classes $\seq{y}{\ell}$ extend to global classes on 
$(G/T)_{T}$.  Indeed, since $y_{i}= c_{1}(L_{\chi_{i}})$ is the first 
Chern class of a $T$-equivariant line bundle on $G/T$, 
by \eqref{e:extension} its global 
extension is
\[
\ty_{i} = c_{1}((L_{\chi_{i}})_{T})
\]
in $H_T^{*}(G/T)$.  Thus, $H_T^{*}(G/T)$ is generated as a $\Q$-algebra by 
$\seq{u}{\ell}$, $\seq{\ty}{\ell}$, and it remains only to determine the 
relations they satisfy.

\section{Restriction to a fixed point in $G/T$}

Although it is possible to give a shorter derivation of the equivariant
cohomology ring $H_T^*(G/T)$ (see \cite{brion}),
we will determine the ring structure of $H_T^*(G/T)$
from the fixed points of the
action of $T$ on $G/T$ using the 
following localization theorem of Borel.
This approach leads to a restriction formula
(Proposition~\ref{p:restrictgt})
that will be useful in our subsequent computation of characteristic numbers.

\medskip
\noindent
{\bf Theorem} (Borel localization theorem, \cite[Proposition 2, p.~39]{hsiang})
{\it Suppose a torus $T$ acts on a manifold $M$ with fixed point set
  $F$. 
Let $i_F\colon F \hookrightarrow M$ be the inclusion map.  
Then both the kernel and cokernel of the restriction 
homomorphism
\[
i_F^*\colon H_T^*(M) \rightarrow H_T^*(F)
\]
are torsion $H^*(BT)$-modules.
}

To apply this theorem to the action of $T$ on $G/T$ by left multiplication,
it is necessary to know the fixed point set $F$ of the action
as well as the restriction homomorphism $i_F^*\colon H_T^*(G/T)
\rightarrow H_T^*(F)$.

\begin{proposition} \label{p:fixed}
For a compact Lie group $G$ with maximal torus $T$,
	let $N_G (T)$ be the normalizer of $T$ in $G$ and let $T$ act on 
	$G/T$ by left multiplication.  Then the fixed point set $F$ 
of the action 
	is $W_G =N_G (T)/T$, the Weyl group of $T$ in $G$.
\end{proposition}

\begin{proof}
The coset	$xT$ is a fixed point
\begin{align*}
	\text{iff }\quad & txT=xT  \quad\text{for all } t\in T \\
	\text{iff }\quad & x\inv txT= T  \quad\text{for all } t\in T \\
	\text{iff }\quad & x\inv tx \in T  \\
	\text{iff }\quad & x \in N_G (T). \tag*{\qed}
\end{align*}
\end{proof}

At a fixed point $w$ in $G/T$, the group $T$ acts on the fiber
$(L_{\gamma})_w$ of the line bundle $L_{\gamma}$.
Therefore, the fiber $(L_{\gamma})_w$ is a complex representation of $T$.

\begin{lemma} \label{l:fiber}
	At the fixed point $w=xT \in W_G$, the torus $T$ acts on the fiber of 
the line bundle 
	$L_{\gamma}$ as the representation $w\cdot \gamma$, i.e.,
$(L_{\gamma})_w = \C_{w\cdot\gamma}$.
\end{lemma}

\begin{proof} 
	The fiber of $L_{\gamma}$ at a fixed point $xT$ consists of elements 
	of the form $[x,v] \in G\times_{\gamma}\C$. If $t\in T$, then
\begin{align*}
	t\cdot [x,v] &= [tx, v]  = [x(x\inv t x), v] \\
	&= [x, \gamma (x\inv t x) v] = [x, (w\cdot \gamma)(t) v].
\end{align*}
Hence, under the identification $[x,v] \leftrightarrow v$, the torus $T$ 
acts on the fiber $(L_{\gamma})_w$ as the representation $w\cdot \gamma$.
\qed\end{proof}

\begin{lem} \label{l:resbundle}
Let $w$ be a point in the Weyl group $W_G \subset G/T$
and $i_w\colon \{ w\} \hookrightarrow G/T$ the inclusion map.
For a character $\gc$ of $T$, the restriction of the line 
bundle $(L_{\gc})_T$ from 
$(G/T)_T$ to $\{w\}_T \simeq BT$ is given by
\[
(i_w)_T^* (L_{\gc})_T \simeq S_{w\cdot \gc},
\]
where $S_{w\cdot \gc}$ is the complex line bundle on $BT$
 associated to the character $w\cdot \gc$.
\end{lem}

\begin{proof}
By Lemma \ref{l:fiber} the restriction of the line bundle $L_{\gc}$ to the fixed 
point $w$ gives rise to a commutative diagram of 
$T$-equivariant maps
\[
\bfig
\Square/^{ (}->`>`>`^{ (}->/<500>%
[\C_{w\cdot \gc}`L_{\gc}`\{w\}`G/T.;```i_w]
\efig
\]
Taking the homotopy quotients results in the
diagram
\[
\bfig
\Square/^{ (}->`>`>`^{ (}->/<500>%
[(\C_{w\cdot \gc})_T`(L_{\gc})_T`BT`(G/T)_T.;```(i_w)_T]
\efig
\]
But $(\C_{w\cdot \gc})_T = ET \times_T \C_{w\cdot \gc}$
is precisely the line bundle $S_{w\cdot \gc}$ over $BT$ 
associated to the character $w\cdot \gc$ of $T$.
\qed\end{proof}

To avoid a plethora of subscripts, we write $(i_w)_T^*$ also as
$i_w^*$.
To describe the restriction $i_F^*\colon H_T^*(G/T) \rightarrow H_T^*(F)$,
we need to describe $i_w^*(u_i)$ and $i_w^*(\tilde{y}_i)$ for each $w
\in W$.

\begin{prop}[Restriction formula for $G/T$] \label{p:restrictgt}
At a fixed point $w \in W$, let 
\[
i_w^*\colon H_T^*(G/T) \rightarrow H_T^*(\{w\})
\]
be the restriction map
in equivariant cohomology. Then
\begin{enumerate}
\item[(i)] $i_w^* u_i = u_i$.
\item[(ii)] For any character $\gamma \in \hat{T}$, 
$i_w^* c_1^T(L_{\gamma}) = w
\cdot c_1(S_{\gamma})$.
In particular, if
$\tilde{y}_i = c_1^T(L_{\chi_i})$ and $u_i = c_1(S_{\chi_i})$,
then
$ i_w^* \tilde{y}_i= w\cdot u_i$.
\end{enumerate}
\end{prop}

\begin{proof}
(i) Let $\pi\colon G/T \rightarrow \{w\}$ 
be the constant map.
Since $\pi\comp i_w = 1_w$, the identity map on $\{w\}$,
$i_w^* \pi^* = 1_{H_T^*(G/T)}$.
The elements $u_i$ in $H_T^*(G/T)$ are really $\pi^*u_i$, so
\[
i_w^* u_i = i_w^* \pi^* u_i =u_i.
\]
(ii) By the naturality of $c_1$ and by
Lemma \ref{l:resbundle},
\begin{align*}
i_w^* c_1((L_{\gamma})_T) &= c_1(i_w^* (L_{\gamma})_T)
=  c_1(S_{w\cdot \gamma}) \\
 &=w\cdot c_1(S_{\gamma}) \qquad \text{(by 
 Corollary~\ref{c:weylaction}).}
 \end{align*}
Now take $\gamma$ to be $\chi_i$.
%Since $i_w^*$ and $w$ are both algebra homomorphisms, for 
%any polynomial $f$ in $\ell$ variables,
%\[
%i_w^* f(y_1, \dots, y_{\ell}) = w\cdot f(u_1, \cdots, 
%u_{\ell}).
%\]
\qed\end{proof}

\section{The equivariant cohomology rings of $G/T$ and $G/H$}

Suppose $G$ is a compact, connected Lie group with maximal torus $T$.
Choose a basis $\chi_1, \ldots, \chi_{\ell}$ of the character group
$\hat{T}$.
Let $L_{\chi_i}$ be the associated line bundles over $G/T$,
and $y_i =c_1(L_{\chi_i}) \in H^2(G/T)$ and $\tilde{y}_i =
c_1^T(L_{\chi_i}) \in H_T^2(G/T)$ be their ordinary and equivariant first
Chern classes.
Similarly, let $S_{\chi_i}$ be the associated line bundles over $BT$,
and $u_i = c_1(S_{\chi_i}) \in H^2(BT)$ their first Chern classes.
Recall $R = H^*(BT) = \Q[u_1, \ldots, u_{\ell}]$.

\begin{thm} \label{t:ring}
\begin{enumerate}
\item[(i)] The equivariant cohomology ring of $G/T$ under the action of $T$ on
$G/T$ by left multiplication is
\[
H_T^*(G/T) \simeq \frac{\Q[u_1, \ldots, u_{\ell}, \tilde{y}_1, \ldots,
\tilde{y}_{\ell}]}{\mathcal{J}},
\]
where $\mathcal{J}$ is the ideal in $\Q[u_1, \ldots, u_{\ell}, \tilde{y}_1,
  \ldots,
\tilde{y}_{\ell}]$ generated by $b(\tilde{y})-b(u)$ for all
polynomials
$b \in R_{+}^{W_G}$, the $W_G$-invariant polynomials of positive degree.
\item[(ii)] If $H$ is a closed subgroup of $G$ containing the maximal torus $T$,
and $T$ acts on $G/H$ by left multiplication,
then
\[
H_T^*(G/H) \simeq \frac{\Q[u_1, \ldots, u_{\ell}] \otimes_{\Q}
(\Q[\tilde{y}_1, \ldots, \tilde{y}_{\ell}]^{W_H})}{\mathcal{J}},
\]
where $\mathcal{J}$ is the ideal in $\Q[u_1, \ldots, u_{\ell}] \otimes_{\Q}
(\Q[\tilde{y}_1, \ldots, \tilde{y}_{\ell}]^{W_H})$
generated by $b(\tilde{y})-b(u)$ for all
polynomials
$b \in R_{+}^{W_G}$.
\end{enumerate}
\end{thm}

\begin{proof}
By Proposition~\ref{p:fixed}, the fixed point set $F$ of the action of $T$
on $G/T$ is the Weyl group $W = W_G$.
For each $w \in W$, by the restriction formula
(Proposition~\ref{p:restrictgt}), $i_w^* \tilde{y}_i = w\cdot u_i$.
Hence, if $b(u)= b(\seq{u}{\ell})$ is a $W$-invariant polynomial with
coefficients in $\Q$, then
\begin{align*}
    i_{w}^{*}b(\seq{\ty}{\ell})&= b(w\cdot u_{1}, \ldots, w\cdot u_{\ell}) \\
    &=b(\seq{u}{\ell}).
\end{align*}
With $\pi\colon G/T \rightarrow \{ w\}$ being the constant map and
$i_w\colon \{ w\} \rightarrow G/T$ the inclusion map,
$\pi \comp i_w = 1$.
Thus,
\[
b(u) = i_w^*\pi^* b(u).
\]
It follows that
\[
i_w^* b(\ty) = b(u) = i_w^*\pi^*b(u),
\]
or
\[
i_w^*(b(\ty) - \pi^*b(u))=0.
\]
As is customary, in $H_T^*(G/T)$, we identify $\pi^*b(u)$ with $b(u)$,
so we can write
\[
i_w^*(b(\ty)-b(u)) = 0.
\]
Since this is true at any fixed point $w \in F$,
\[
i_{F}^{*} (b(\ty)-b(u))=0,
\]
where $i_{F}\colon F \hookrightarrow G/T$ is the inclusion map.

Let $R = H^*(BT)$.
By the Borel localization theorem, the kernel of the restriction map
\[
i_{F}^{*}: H_{T}^{*} (G/T) \rightarrow H_{T}^{*}(F)
\]
is a torsion $R$-module.  Since $H_{T}^{*}(G/T)$ is a free $R$-module, 
the kernel must be $0$.  
Therefore, $i_{F}^{*}$ is injective.
So we obtain the relations
\begin{equation} \label{e:relations}
b(\ty) - b(u) =0
\end{equation}
in $H_{T}^{*}(G/T)$ for all $W$-invariant polynomials $b(\seq{u}{\ell})$.

Let $\mathcal{J}$ be the ideal $(b(\ty) - b(u))$ 
in $\Q[\seq{u}{\ell}, \seq{\ty}{\ell}]$ generated by 
$b(\ty) - b(u)$ for all
polynomials $b\in R_+^W$.  Then there is a 
surjective ring homomorphism
\begin{equation} 
   \phi\colon  \dfrac{\Q[\seq{u}{\ell}, \seq{\ty}{\ell}]}{\mathcal{J}}
\simeq \frac{R \otimes_{\Q} R}{\mathcal{J}} \rightarrow H_{T}^{*}(G/T).
\end{equation}

From the spectral sequence of the fibering 
\eqref{e:fibering},  we know that $H_T^*(G/T)$ is a free $R$-module of rank 
equal to $\dim H^*(G/T)$.  
To prove that $\phi$ is an isomorphism, it suffices to show that
$(R\otimes_{\Q} R)/\mathcal{J}$ is a free $R$-module of the same rank,
for in that case the kernel of $\phi$, being of rank 0, is a torsion
submodule
of a free module and is therefore the zero module.
Note that
\[
\dfrac{R\otimes_{\Q}R}{\mathcal{J}} \simeq R \otimes_{R^W} R,
\]
since the tensor product over $R^W$ is obtained from the tensor product over $\Q$ 
by dividing out by the ideal generated by elements of the form $b\otimes 
1 - 1 \otimes b$ for $b \in R^W$.
Moreover,
\begin{align*}
R \otimes_{R^W} R &\simeq (R\otimes_{\Q} \Q)\otimes_{R^W} R 
\simeq \left( R \otimes_{\Q} \dfrac{R^W}{(R_+^W)}\right) \otimes _{R^W} R \\
&\simeq R \otimes_{\Q} \left( \dfrac{R^W}{(R_+^W)} \otimes_{R^W} R\right) 
\simeq R \otimes_{\Q} \dfrac{R}{(R_+^W)} \\
&\simeq R \otimes_{\Q} H^*(G/T) \qquad \text{(Theorem \ref{t:cohomgt}).}
\end{align*}
So $(R\otimes_{\Q} R)/\mathcal{J}$ is a free $R$-module of rank equal to $\dim 
H^*(G/T)$. 
This proves that $\phi$ is an isomorphism.  

(ii) Since the fiber bundle $G/H \rightarrow G/T$ with fiber $H/T$
is $T$-equivariant,
it induces a map in homotopy quotients $(G/H)_T
\rightarrow (G/T)_T$, which is also a fiber bundle with fiber $H/T$.
By Lemma \ref{l:gmodtbundle},
\begin{align} %\label{e:equivcohomgh}
H_T^*(G/H) &= H_T^*(G/T)^{W_H} \notag\\
&= \dfrac{\Q [\seq{u}{\ell}] \otimes_{\Q} (\Q [\seq{\ty}{\ell}]^{W_H})}{\mathcal{J}},
\end{align}
where as before, $\mathcal{J}$ is the ideal generated by $b(\ty)-b(u)$ for all $b \in 
R_+^{W_G}$.
\qed\end{proof}

\section{The tangent bundle of $G/H$} \label{s:tangent}

Let $G$ be a Lie group with Lie algebra $\frakg$, and $H$ a closed 
subgroup with Lie algebra $\frakh$.  

% If $\pi\colon E \rightarrow G/H$ is a $G$-equivariant vector bundle, then 
% at the neutral element $H \in G/H$, the fiber $E_0 := 
% \pi\inv(\{H\})$ is a representation of $H$.  In 
% (\cite{segal}, p.~130) Segal shows that $E$ is isomorphic to 
% the associated bundle $G \times_H E_0$.  
% In particular, $T(G/H) \simeq G \times_H (\frakg/\frakh)$, 
% where $H$ acts on $\frakg/\frakh$ by left multiplication, or 
% what is the same in this case, the adjoint representation.
% We give a direct proof below.

\begin{prop}[See also {\cite[p.~130]{segal}}] \label{p:tangh}
The tangent bundle of $G/H$ is diffeomorphic to the vector bundle 
$G\times_H (\frakg/\frakh)$ associated to the principal $H$-bundle $G \rightarrow 
G/H$ via the adjoint representation of $H$ on $\frakg/\frakh$.
\end{prop}

\begin{proof}
For $g\in G$, let ${\ell}_g\colon G/H \rightarrow G/H$ denote 
left multiplication by 
$g$.  Then ${\ell}_g$ is a diffeomorphism and its differential
\[
({\ell}_{g})_*\colon T_{eH}(G/H) \rightarrow T_{gH}(G/H)
\]
is an isomorphism.  

Now define
\[
\varphi\colon G \times _H (\frakg/\frakh) \rightarrow T(G/H)
\]
by
\[
\varphi([g,v]) = ({\ell}_g)_* (v).
\]
To show that $\varphi$ is well-defined, pick another representative of 
$[g,v]$, say $[gh, \Ad(h\inv)v]$ for some $h \in H$.
Then
\begin{align*}
({\ell}_{gh})_*(\Ad(h\inv)v) &= ({\ell}_{gh})_* ({\ell}_{h\inv})_* (r_h)_*v\\
&= ({\ell}_g)_*({\ell}_{h})_*({\ell}_{h\inv})_* (r_h)_*v\\
&=({\ell}_g)_*v,
\end{align*}
since right multiplication $r_h$ on $G/H$ is the identity map.

Since $\dim G/H = \dim \frakg/\frakh$,
$\varphi$ is a surjective bundle map of two 
vector bundles of the same rank, and so
it is an isomorphism.
\qed\end{proof}

\section{Equivariant Euler class of the normal bundle at a 
fixed point of $G/T$}

In this section $G$ is a compact, connected Lie group, $T$ a 
maximal torus in $G$, 
$\frakg$ and $\frakt$ their 
respective Lie algebras,
and $W= N_G(T)/T$ the Weyl group of $T$ in $G$.  
The adjoint representation of $T$ 
on $\frakg$ decomposes $\frakg$ into a direct sum
\[
\frakg = \frakt \oplus \bigg(\bigoplus_{\ga \in \tri^+} \C_{\ga}\bigg),
\]
where $\tri^+$ is a choice of positive roots.

Let $w\in W$ be a fixed point of the left multiplication action 
of $T$ on $G/T$, and $\nu_w$ the normal bundle of $\{w\}$ 
in $G/T$.  The normal bundle at a point is simply the 
tangent space. Let $\ell(w)$ be the length of $w$ and 
$(-1)^w := (-1)^{\ell(w)}$ the sign of $w$.

\begin{prop}[Euler class formula for $G/T$] \label{p:eulergt}
The equivariant Euler class of the normal bundle $\nu_w$ at 
a fixed point $w \in W$ of the left action of $T$ on $G/T$ is
\[
e^T(\nu_w) = w\cdot \left( \prod_{\ga\in \tri^+} 
c_1(S_{\ga})\right) 
= (-1)^w \prod_{\ga \in \tri^+} c_1(S_{\ga}) 
\in H^*(BT).
\]
\end{prop}

\begin{proof}
By Proposition~\ref{p:tangh} the tangent bundle of $G/T$ 
is the homogeneous vector bundle associated to the adjoint 
representation of $T$ on $\frakg/\frakt = 
\oplus_{\ga \in \tri^+} \C_{\ga}$.
In the notation of Section \ref{s:assobundles},
with $L_{\alpha}$ being the complex line bundle over $G/T$
associated to the character $\alpha$ of $T$,
\begin{equation} \label{e:tangt}
T(G/T) \simeq G \times_T (\frakg/\frakt) \simeq
 G \times_T \left(\bigoplus_{\ga \in \tri^+} 
\C_{\ga}\right) \simeq \bigoplus_{\ga \in \tri^+} L_{\ga}.
\end{equation}

By \eqref{e:tangt} and
Lemma \ref{l:fiber}, at a fixed point $w$ the normal bundle 
$\nu_w$ is
\[
\nu_w = T_w(G/T) = \bigoplus_{\ga \in \tri^+} (L_{\ga})_w \simeq
\bigoplus_{\ga \in \tri^+} \C_{w\cdot \ga}.
\]
It follows that the equivariant Euler class of $\nu_w$ is
\begin{alignat*}{2}
e^T(\nu_w) &= e^T\left(\bigoplus_{\ga \in \tri^+} \C_{w\cdot\ga}\right) =
e\left(\bigoplus_{\ga \in \tri^+} (\C_{w\cdot\ga})_T\right) \\
&= e\left(\bigoplus_{\ga \in \tri^+} S_{w\cdot \ga}\right) 
&&\qquad\text{(by the definition of $S_{\ga}$)}\\
&=\prod_{\ga \in \tri^+} c_1(S_{w\cdot \ga}) = 
\prod_{\ga \in \tri^+} w\cdot c_1(S_{\ga}) 
&&\qquad\text{(Corollary~\ref{c:weylaction})}\\
&=w\cdot \left(\prod_{\ga \in \tri^+}c_1(S_{\ga})\right)
&&\qquad\text{($w$ is 
a ring homomorphism)}
\end{alignat*}

Each simple reflection $s$ in the Weyl group $W=N_G(T)/T$ 
carries exactly one 
positive root to a negative root.  
Note that $c_1(S_{-\alpha}) = c_1(S_{\alpha}^{\spcheck}) =
-c_1(S_{\alpha})$.
Hence,
\[
s\cdot \left(
\prod_{\ga \in \tri^+} c_1(S_{\ga}) \right) =
\prod_{\ga \in \tri^+} c_1(S_{s\cdot \ga}) 
=-\prod_{\ga \in \tri^+} c_1(S_{\ga}).
\]
Since $w$ is the product of $\ell (w)$ reflections, 
\[
w\cdot 
\left(
\prod_{\ga \in \tri^+} c_1(S_{\ga}) \right) 
=(-1)^{\ell(w)}\prod_{\ga \in \tri^+} c_1(S_{\ga}). \tag*{\qed}
\]
\end{proof}

\section{Fixed points of $T$ acting on $G/H$}

The quotient space $W_G /W_H$ can be viewed as a 
subset of $G/H$ via
\[
\frac{W_G}{W_H} \simeq \frac{N_G (T)}{N_H (T)} = \frac{N_G(T)}{N_G(T) \cap 
H} 
\hookrightarrow \frac{G}{H}.
\]

\begin{prop} \label{p:fixedgh}
Under the action of $T$ on $G/H$ by left multiplication, the fixed point 
set $F$ is $W_G/W_H$.
\end{prop}

\begin{proof}
The coset $xH$ is a fixed point
\begin{align*}
	\text{iff }\quad & txH=xH  \qquad\text{for all } t\in T \\
	\text{iff }\quad & x\inv txH = H  \qquad\text{for all } t\in T \\
	\text{iff }\quad & x\inv tx \in H  \qquad\text{for all } t\in T \\
	\text{iff }\quad & x\inv Tx \subset H .
\end{align*}
Since any two maximal tori in $H$ are conjugate by an element of $H$, 
there is an element $h \in H$ such that
\[
h\inv x\inv T xh = T.
\]
Therefore, $xh \in N_G (T)$, and $xH=xhH$.  Thus any fixed point can be 
represented as $yH$ for some $y \in N_G(T)$.

Conversely, if $y \in N_G(T)$, then $yH$ is a fixed point 
of the action of $T$ on $G/H$, since $TyH=yTH=yH$.  
It follows that there 
is a surjective map
\begin{align*}
N_G(T) &\twoheadrightarrow F \subset G/H, \\
y& \mapsto yH
\end{align*}
with fiber $N_G(T) \cap H$ 
and hence a bijection
\[
\dfrac{N_G(T)}{N_G(T)\cap H} \simeq F. \tag*{\qed}
\]
\end{proof}

\section{Restriction to a fixed point of $G/H$}

\begin{thm}[Restriction formula for $G/H$] \label{t:restrictgh}
With $H_T^*(G/H)$ as in Theorem \ref{t:ring}(ii) and 
$w$ a fixed point in $W_G/W_H$, 
the restriction map $i_w^*$ in equivariant cohomology
\[
i_w^*\colon H_T^*(G/H) \rightarrow H_T^*(\{w\}) = \Q[u_1, \dots, 
u_{\ell}]
\]
is given by
\[
i_w^*u_i = u_i, \qquad i_w^*f(\ty) = w\cdot f(u)
\]
for any $W_H$-invariant polynomial $f(\ty) \in \Q[\ty_1, \dots, 
\ty_{\ell}]^{W_H}$.
\end{thm}

\begin{proof}
Write $w=xH \in N_G(T)/(N_G(T) \cap H) \subset G/H$, with $x \in N_G(T)$.
The commutative diagram of $T$-equivariant maps 
\[
\bfig
\Square/^{ (}->`>`>`^{ (}->/%
[\{xT\}`G/T`\{xH\}`G/H;{i_{xT}}`{\gs}_x`\gs`{i_{xH}}]
\efig
\] 
induces a commutative diagram in equivariant cohomology
\[
\bfig
\Square/<-`<-`<-`<-/%
[H^*(BT)`H_T^*(G/T)`H^*(BT)`H_T^*(G/H)\ .;i_{xT}^*`{\gs}_x^*`\gs^*`i_{xH}^*]
\efig
\]
Since $\gs^*\colon H_T^*(G/H) \rightarrow H_T^*(G/T)$ is an injection
and 
\[
\sigma_x^*\colon H_T^*(\{xH\}) =H^*(BT) \rightarrow
H^*(BT)=H_T^*(\{xT\})
\] 
is the identity map,
the restriction $i_{xH}^*$ is given by the same formula as 
the restriction $i_{xT}^*$.  The theorem then follows from the
restriction formula for $G/T$ (Proposition~\ref{p:restrictgt}).
\qed\end{proof}

\section{Pullback of an associated bundle}

Assume that $\rho\colon H \rightarrow \GL(V)$ is a representation of
the group $H$ on a vector space $V$ over any field.
By restriction, one obtains a representation of the maximal 
torus $T$.

\begin{prop} \label{p:pullbackass}
Let $\gs\colon G/T \rightarrow G/H$ be the projection map.  Under $\gs$ 
the associated bundle $G\times_H V$ pulls back to $G\times_T 
V$:
\[
\gs^*(G\times_H V) \simeq G\times_T V.
\]
\end{prop}

\begin{proof}
Let $[g,v]_T$ and $[g,v]_H$ denote the equivalence classes 
of $(g,v)$ in $G\times_TV$ and $G\times_H V$ respectively.
Define $\varphi\colon G\times_T V \rightarrow \gs^*( G\times_H V)$ by
\[
[g,v]_T \mapsto (gT, [g,v]_H).
\]

If $(gT, [g',v']_H)$ is any element of $\sigma^*(G\times_H 
V)$, then $gH= g'H$.  Hence, $g'=gh$ for some $h \in H$ and
\begin{align*}
(gT, [g',v']_H) &= (gT, [gh,v']_H) \\
&=(gT, [g, hv']_H) \\
&= \varphi ([g, hv']_T),
\end{align*}
which shows that $\varphi$ is surjective.  A surjective 
bundle map between two vector bundles of the same rank is an 
isomorphism.
\qed\end{proof}

\section{Pulling back the tangent bundle of $G/H$ to $G/T$}

Let $\frakg, \frakh$, and $\frakt$ be the Lie algebras of 
$G, H$, and $T$ respectively.  Under the adjoint 
representation of $H$ the Lie algebra $\frakg$ decomposes 
into a direct sum of $H$-modules
\[
\frakg = \frakh \oplus \frakm.
\]
By Proposition~\ref{p:tangh}, the tangent bundle of $G/H$ 
is the associated bundle
\begin{equation} \label{e:tangentass}
T(G/H) \simeq G\times_H \frakm.
\end{equation}

Restricting the adjoint representation to the maximal 
torus $T$, the Lie algebras of $H$ and $G$ decompose further 
into a sum of $T$-modules
\[
\frakh = \frakt \oplus \left(\bigoplus_{\ga \in \tri^+(H)} \C_{\ga}\right)
\]
and
\[
\frakg = \frakt \oplus \left(\bigoplus_{\ga \in \tri^+(H)} \C_{\ga}\right)
\oplus \left(\bigoplus_{\ga \in \tri^+ \setminus \tri^+(H)} \C_{\ga}\right),
\]
where $\tri^+(H)$ denotes a choice of positive roots for 
$H$, $\tri^+$ a choice of positive roots for $G$ 
containing $\tri^+(H)$, and $\tri^+ \setminus \tri^+(H)$ 
the complement of $\tri^+(H)$ in $\tri^+$.  Hence, as a $T$-module,
\begin{equation} \label{e:m}
\frakm = \bigoplus_{\ga \in \tri^+ \setminus \tri^+(H)} \C_{\ga}.
\end{equation}

\begin{prop} \label{p:pullbacktangh}
Under the natural projection $\gs\colon G/T \rightarrow G/H$, the 
tangent bundle $T(G/H)$ pulls back to a sum of associated 
line bundles:
\[
\gs^*T(G/H) \simeq \bigoplus_{\ga \in \tri^+ \setminus \tri^+(H)} 
L_{\ga}.
\]
\end{prop}

\begin{proof}
\begin{alignat*}{2}
\gs^*T(G/H) &\simeq \gs^*(G \times_H \frakm) \simeq G 
\times_T \frakm &&\qquad\text{(by \eqref{e:tangentass} and 
Proposition~\ref{p:pullbackass})}\\
&= G\times_T \left(\bigoplus_{\ga \in \tri^+ \setminus \tri^+(H)}
\C_{\ga}\right)
&&\qquad\text{(by \eqref{e:m})}\\
&= \bigoplus_{\ga \in \tri^+ \setminus \tri^+(H)} L_{\ga}&&
\qquad\text{(definition of $L_{\ga}$)} \tag*{\qed}
\end{alignat*}
\end{proof}

\section{The normal bundle at a fixed point of $G/H$}

At a fixed point $w=xH$ of $G/H$, with $x\in N_G(T)$, the 
action of $T$ on $G/H$ induces an action of $T$ on the 
tangent space $T_w(G/H)$.

\begin{prop} \label{p:tanspgh}
At a fixed point $w=xH$ of $G/H$, the tangent space 
$T_w(G/H)$ decomposes as a representation of $T$ into
\[
T_w(G/H) \simeq \bigoplus_{\ga \in \tri^+ \setminus \tri^+(H)} 
\C_{w\cdot \ga}.
\]
\end{prop}

\begin{proof}
Let $\gs\colon G/T \rightarrow G/H$ be the projection map.  Then
\begin{alignat*}{2}
T_{xH}(G/H) &= (T(G/H))_{xH} &&\qquad\text{(fiber of 
tangent bundle at $xH$)}\\
&=(\gs^*T(G/H))_{xT}&&\\
&=\left(\bigoplus_{\ga \in \tri^+ \setminus \tri^+(H)} (L_{\ga})_{xT}\right)
&&\qquad\text{(Proposition~\ref{p:pullbacktangh})}\\
&\simeq \bigoplus_{\ga \in \tri^+ \setminus \tri^+(H)} \C_{w\cdot\ga}
&&\qquad\text{(Lemma \ref{l:fiber})}  \tag*{\qed}
\end{alignat*}
\end{proof}

\begin{thm}[Euler class formula for $G/H$] \label{t:eulergh}
At a fixed point $w \in W_G/W_H$ of the left action of $T$ 
on $G/H$, the equivariant Euler class of the normal bundle 
$\nu_w$ is given by
\[
e^T(\nu_w) = w\cdot\left( \prod_{\ga \in \tri^+ \setminus \tri^+(H)} 
c_1(S_{\ga}) \right).
\]
\end{thm}

\begin{proof}
The normal bundle $\nu_w$ at the point $w$ is the tangent 
space $T_w(G/H)$.  
By the multiplicativity of the Euler class and the fact that the Euler
class of a complex line bundle is its first Chern class,
the equivariant Euler class of $\nu_w$ is
\begin{align*}
e^T(\nu_w) &= e^T (T_w(G/H)) = e^T \left(\bigoplus_{\ga \in \tri^+ \setminus \tri^+(H)}
\C_{w\cdot \ga}\right) \qquad \text{(Proposition~\ref{p:tanspgh})} \\
&= \prod_{\ga \in \tri^+ \setminus \tri^+(H)} c_1^T(\C_{w\cdot 
\ga}) = \prod_{\ga \in \tri^+ \setminus \tri^+(H)} 
 c_1((\C_{w\cdot\ga})_T) \\
&=\prod_{\ga \in \tri^+ \setminus \tri^+(H)} c_1(S_{w\cdot\ga}) = w\cdot \left(
\prod_{\ga \in \tri^+ \setminus \tri^+(H)} 
c_1(S_{\ga})\right).  \tag*{\qed}
\end{align*}
\end{proof}

\section{Equivariant characteristic numbers of $G/H$ and $G/T$} 
\label{s:eqcharnogh}

Suppose a torus $T$ acts on a compact oriented manifold $M$.
Let $\pi\colon M \rightarrow \pt$ be the constant map and 
$\pi_*\colon H_T^*(M;\R) \rightarrow H_T^*(\pt;\R)= H^*(BT;\R)$ the push-forward 
or integration
map in equivariant cohomology.
For any class $\tilde{\eta} \in H_T^*(M;\R)$, the equivariant
 localization formula computes 
$\pi_*\tilde{\eta}$ in terms of an integral over the fixed point set 
$F$ (\cite{atiyah-bott}, \cite{berline-vergne}).  In case 
the fixed points are isolated, the equivariant localization formula 
states that
\begin{equation} \label{e:localization}
\pi_*\tilde{\eta} = \sum_{p \in F} \dfrac{i_p^*\tilde{\eta}}{e^T(\nu_p)}.
\end{equation}
On the right-hand side the calculation is performed in the fraction
field of $H^*(BT;\R)$ and so is a priori a rational function of
$u_1, \ldots, u_{\ell}$, but it is part of the theorem that the sum
will be a polynomial in $u_1, \ldots, u_{\ell}$ since the left-hand
side $\pi_*\tilde{\eta} \in H^*(BT;\R)$ is a polynomial in $u_1, \ldots,
u_{\ell}$. 

Although the equivariant localization formula
is stated for real cohomology, by viewing
rational cohomology as a subset of real cohomology, we can apply the
equivariant localization formula to rational cohomology classes.
Now suppose $G$ is a compact, connected Lie group, $T$ a maximal torus in $G$, and $H$ a closed subgroup of $G$ containing $T$.
For $\tilde{\eta} \in H_T^*(G/H;\Q)$,
$\pi_* \tilde{\eta}$ is a priori a real class in $H^*(BT;\R)$.
However, the explicit formula \eqref{e:gysingh} below
shows that it is in fact
a rational class in $H^*(BT;\Q)$.

For $M=G/H$ with the standard torus action, 
by Theorem \ref{t:ring}(ii), an element of $H_T^*(G/H;\Q)$ is of the form
\[
\tilde{\eta} = \sum (\pi^*a_i(u)) f_i(\ty),
\]
where $a_i(u) \in \Q[u_1, \ldots, u_{\ell}]$ and
$f_i(\ty) \in \Q[\ty_1, \ldots, \ty_{\ell}]^{W_H}$.
By the projection formula,
\[
\pi_* \tilde{\eta} = \sum a_i(u) \pi_* f_i(\ty).
\]
Thus, to calculate $\pi_*\colon H_T^*(G/H;\Q) \subset H_T^*(G/H;\R)
\rightarrow H^*(BT;\R)$, it suffices to calculate
$\pi_* f(\ty)$ for $f(\ty)= f(\ty_1, \ldots, \ty_{\ell})$
a $W_H$-invariant polynomial with coefficients in $\Q$.
Since $\ty_1, \ldots, \ty_{\ell}$ are all equivariant 
characteristic classes, 
$\pi_*f(\ty)$ is called an \emph{equivariant characteristic 
number} of $G/H$.
With the aid of 
 the restriction formula (Theorem~\ref{t:restrictgh}) and 
the Euler class formula (Theorem~\ref{t:eulergh}),
the equivariant localization formula
\eqref{e:localization} gives
for any $f(\ty) \in \Q[\ty_1, \ldots, \ty_{\ell}]^{W_H}$,
\begin{equation} \label{e:gysingh}
\pi_* f(\ty) = \sum_{w\in W_G/W_H} w \cdot \left(
\dfrac{f(u)}{\prod_{\ga \in \tri^+ \setminus \tri^+(H)} c_1(S_{\ga})}
\right).
\end{equation}
In this formula, 
$f(u) = f(u_1, \ldots, u_{\ell})$ is obtained from $f(\tilde{y})$
by replacing $\tilde{y}_i$ by $u_i$.
Since the left-hand side $\pi_*f(\ty) \in H^*(BT;\R)$ is a polynomial in
$u_1, \ldots, u_{\ell}$ with real coefficients, so is the right-hand
side.
But the right-hand side clearly has rational coefficients.
Hence, $\pi_*f(\ty) \in H^*(BT;\Q)$.

For $G/T$, 
\[
H_T^*(G/T) = \dfrac{\Q[u_1, \dots, u_{\ell}, \ty_1, \dots, 
\ty_{\ell}]}
{(b(\ty)-b(u) \mid b \in R_+^{W_G})},
\]
and the fixed point set is $W_G$.
For $f(\ty) \in \Q[\ty_1, \ldots, \ty_{\ell}]$,
by the restriction formula (Proposition~\ref{p:restrictgt})
and the Euler class formula (Proposition~\ref{p:eulergt})
for $G/T$,
\begin{equation} \label{e:gysingt}
\pi_* f(\ty) = \sum_{w\in W_G} w \cdot \left(
\dfrac{f(u)}{\prod_{\ga \in \tri^+} c_1(S_{\ga})}
\right)
= \dfrac{\sum_{w\in W_G} (-1)^w w\cdot f(u)}
{\prod_{\ga \in \tri^+} c_1(S_{\ga})}.
\end{equation}

\section{Ordinary integration and equivariant integration}

We state a general principle (Proposition~\ref{p:ordint}), well 
known to the experts, relating ordinary integration and 
equivariant integration. 

\begin{prop} \label{p:ordint}
Let $M$ be a compact oriented manifold of dimension $n$ on 
which a compact, connected Lie group $G$ acts
and let $\pi_*\colon H_G^*(M;\R) \rightarrow H^*(BG;\R)$ be equivariant 
integration.  Suppose a cohomology class $\eta \in 
H^n(M;\R)$ has an equivariant extension $\tilde{\eta} \in 
H_G^n(M;\R)$.  Then
\begin{equation} \label{e:int}
\int_M \eta = \pi_* \tilde{\eta}.
\end{equation}
\end{prop}

\begin{rem}
For a torus action the right-hand side $\pi_*\tilde{\eta}$ 
of \eqref{e:int} can be 
computed using the equivariant localization formula in terms of the 
fixed point set $F$ of $T$ on $M$.  In case the fixed 
points are isolated, this gives
\[
\int_M \eta = \sum_{p\in F} \dfrac{i_p^* \tilde{\eta}}
{e^T (\nu_p)}.
\]
\end{rem}

\medskip
\noindent
{\sc Proof of Proposition~\ref{p:ordint}}.
The commutative diagram
\[
\bfig
\Square/^{ (}->`>`>`^{ (}->/[M`M_G`\pt`BG;i`\tau`\pi`j]
\efig
\]
induces by the push-pull formula (\cite[p.~158]{guillemin-sternberg}) a commutative diagram in 
cohomology
\begin{equation} \label{e:restrictfiber}
\bfig
\Square/<-`>`>`<-/%
[{H^n(M;\R)}`{H_G^n(M;\R)}`{H^0(\pt;\R)}`{H^0(BG;\R).};i^*`\tau_*`\pi_*`j^*]
\efig
\end{equation}
In degree 0, the restriction $j^*\colon H^0(BG;\R) \rightarrow
H^0(\pt;\R)=\R$ is  
an isomorphism.  Hence, if $\eta$ has degree $n$ in 
$H^n(M;\R)$, then by the commutative diagram 
\eqref{e:restrictfiber} and the push-pull formula
\[
\int_M \eta = \tau_* \eta = \tau_* i^* \tilde{\eta} = 
j^* \pi_* \tilde{\eta} = \pi_* 
\tilde{\eta}.   \tag*{\qed}
\]
{\vskip 2ex\par}

Since all ordinary characteristic classes of $G$-vector 
bundles have equivariant extensions, all ordinary 
characteristic numbers of $G$-vector bundles can be 
computed from the equivariant localization formula.

Fix a basis $\chi_1, \dots, \chi_{\ell}$ of the character 
group $\hat{T}$ of the maximal torus $T$ in the compact, connected 
Lie group $G$.  Let $H$ be a closed subgroup containing $T$ 
in $G$.  On $G/T$, we have associated bundles 
$L_{\chi_i} := G \times_T \C_{\chi_i}$.  Let $y_i = 
c_1(L_{\chi_i}) \in H^2(G/T;\Q)$.
Let $R$ be the polynomial ring $\Q[y_1, \dots, y_{\ell}]$.
By Theorem \ref{t:cohomgh},
\[
H^*(G/H;\Q) = \dfrac{R^{W_H}}{(R_+^{W_G})} = \dfrac{\Q[y_1, \dots, 
y_{\ell}]^{W_H}}{(R_+^{W_G})}.
\]

\begin{thm} \label{t:charnogh}
Let $f(y) \in \Q[y_1, \dots, y_{\ell}]^{W_H}$ be a 
$W_H$-invariant polynomial of degree $\dim G/H$, where each 
$y_i$ has degree 2. Then the characteristic number $\int_{G/H} f(y)$ 
of $G/H$ is given by
\[
\int_{G/H} f(y) = \pi_* f(\ty) = \sum_{w\in W_G/W_H}
w\cdot \left(
\dfrac{f(u)}{\prod_{\ga \in \tri^+ \setminus \tri^+(H)} c_1(S_{\ga})}
\right).
\]
\end{thm}

\begin{proof}
Since ${\ty}_i= c_1^T(L_{\chi_i})$ is an equivariant extension 
of $y_i$, the cohomology class $f(\ty) \in H_T^*(G/H;\Q)$ is an
equivariant extension of $f(y)$.
Combining Proposition~\ref{p:ordint} and 
\eqref{e:gysingh}, the formula for the ordinary 
characteristic numbers of $G/H$ follows.
As noted earlier, \eqref{e:gysingh} shows that if $f(\ty)$ is a
rational class, then so is $\pi_*f(\ty)$.
\qed\end{proof}

\section{Example: the complex Grassmannian}

In this example, we work out the $T$-equivariant cohomology ring as
well as 
the characteristic numbers of 
the complex Grassmannian $G(k,n)$ of $k$-planes in $\C^n$.  
As a homogeneous space, $G(k,n)$ can be represented as $G/H$, where $G$ is the 
unitary group $U(n)$ and $H$ is the closed subgroup $U(k) 
\times U(n-k)$.

A maximal torus contained in $H$ is
\[
T= \underbrace{U(1) \times \dots \times U(1)}_n =
\left\{
t=\left.
\begin{bmatrix}
t_1 & & \\
& \ddots & \\
&& t_n
\end{bmatrix}
\ 
\right|\ 
t_i \in U(1)
\right\}.
\]
A basis for the characters of $T$ is $\chi_1, \dots, 
\chi_n$, with $\chi_i (t) = t_i$.
The characters $\chi_i$ define line bundles $S_{\chi_i}$
over the classifying space $BT$.
We let $u_i = c_1(S_{\chi_i}) \in H^2(BT)$.
A choice of positive roots for $G$ and for $H$ is
\begin{align*}
\tri^+ &= \{ \chi_i - \chi_j \mid 1 \le i < j \le n \},\\
\tri^+(H) &=  \{ \chi_i - \chi_j \mid 1 \le i < j \le k \} 
\cup
 \{ \chi_i - \chi_j \mid k+1 \le i < j \le n \}.
\end{align*}
Therefore,
\[
\tri^+ \setminus \tri^+(H) =
 \{ \chi_i - \chi_j \mid 1 \le i \le k, \ k+1 \le j \le n \}.
\]
If $\ga = \chi_i - \chi_j$, then
\begin{equation} \label{e:chernclass}
c_1(S_{\ga}) = c_1 (S_{\chi_i} \otimes \dual{S_{\chi_j}}) = 
c_1(S_{\chi_i}) - c_1(S_{\chi_j}) = u_i - u_j.
\end{equation}

The Weyl groups of $T$ in $G$ and $H$ are
\begin{align*}
W_G & = S_n, \ \text{the symmetric group on $n$ letters},\\
W_H & = S_k \times S_{n-k}.
\end{align*}
(Notation: $S_{\ga}, S_{\chi}$ are line bundles over
$BT$ associated to the characters $\ga$ and $\chi$,
but $S_k$, $S_n$ are symmetric groups.)
A permutation in the symmetric group $S_n$ 
is a bijection $w: \{ 1, \ldots, n\} \rightarrow \{ 1, \ldots, n\}$, 
\begin{equation} \label{e:permutation}
w(1) = i_1, \dots, w(k) = i_k, w(k+1)= j_1, \dots, 
w(n)=j_{n-k},
\end{equation}
where $I=(i_1, \ldots, i_k)$ and $J=(j_1, \ldots, j_{n-k})$ are two \emph{complementary} multi-indices, i.e., $I \cup J = \{ 1, \ldots, n\}$. The equivalence class of $w$ in $S_n/(S_k \times S_{n-k})$ has a unique representative with $I = (i_1 < \cdots < i_k)$ and $J=(j_1 < \cdots < j_{n-k})$ both strictly increasing.

By Theorem \ref{t:ring}(ii), the rational equivariant cohomology ring of the 
Grassmannian $G(k, n)$ under the torus action is
\begin{equation} 
H_T^*(G(k,n))=\dfrac{
\Q [\seq{u}{n}] \otimes_{\Q} (\Q [\seq{\ty}{k},\ty_{k+1}, \dots,
\ty_n]^{S_k \times S_{n-k}})
}{\mathcal{J}},
%{(\prod (1+y_i) - \prod(1+u_i))}
\end{equation}
where $\mathcal{J}$ is the ideal in 
$\Q [\seq{u}{n}] \otimes_{\Q} (\Q [\seq{\tilde{y}}{k},\tilde{y}_{k+1}, \dots,
\tilde{y}_n]^{S_k \times S_{n-k}})$
generated by $b(\ty)-b(u)$ for all symmetric polynomials $b(\ty) \in
\Q[\ty_1, \ldots, \ty_n]$.
Let $\sigma_r$ be the $r$th elementary symmetric polynomial.
Since every symmetric polynomial is a polynomial in the elementary
symmetric polynomials,
$\mathcal{J}$ is also the ideal generated by $\sigma_r(\ty) -\sigma_r(u)$
for $r=1, \ldots, n$.
Thus, we may write
\begin{equation} \label{e:htgkn}
H_T^*(G(k,n))=\dfrac{
\Q [\seq{u}{n}] \otimes_{\Q} (\Q [\seq{\ty}{k},\ty_{k+1}, \dots,
\ty_n]^{S_k \times S_{n-k}})}
{(\prod (1+ \ty_i) - \prod(1+u_i))}
\end{equation}
In this formula, the notation $(p(u,\ty))$ means the ideal generated by the 
homogeneous terms of the polynomial $p(u,\ty)$.

If $S$ and $Q$ are the universal sub- and quotient 
bundles over $G(k,n)$, then setting
\begin{align*}
\ts_r = c_r^T(S) &= \gs_r(\ty_1, \dots, \ty_k), \quad\quad 
\tq_r= c_r^T(Q) = \gs_r(\ty_{k+1}, \dots, \ty_n),\\
\intertext{we have}
1+\ts_1+\dots +\ts_k &= \prod_{i=1}^k (1+ \ty_i),
\quad\text{and}\quad 
1+\tq_1+\dots +\tq_{n-k} = \prod_{i=k+1}^n (1+\ty_i).
\end{align*}
Thus \eqref{e:htgkn} can be rewritten in the 
form
\[
H_T^*(G(k,n))=\dfrac{\Q [\seq{u}{n}, \seq{\ts}{k},\seq{\tq}{n-k}]}
{((1+\ts_1+\dots +\ts_k)(1+\tq_1+\dots +\tq_{n-k}) - \prod(1+u_i))}.
\]
Using the relation
\begin{equation} \label{e:relationq}
1+\sum_{i=1}^{n-k} \tq_i = \dfrac{\prod (1+u_i)}{1+\sum \ts_i},
\end{equation}
one can eliminate all the $\tq_i$ from $H_T^*(G(k,n))$; in other words, 
$H_T^*(G(k,n))$ is generated as an algebra over $\Q [\seq{u}{n}]$ by 
$\seq{\ts}{k}$ with relations given by terms of degree $> 2(n-k)$ in 
\eqref{e:relationq}.  In computing degrees, keep in mind that $\deg u_i =2$ 
and $\deg \ts_i = \deg \tq_i = 2i$.

Similarly, in this notation, the rational cohomology ring of the
Grassmannian $G(k,n)$ is
\begin{align*}
H^*(G(k,n)) &=\dfrac{\Q[y_1, \ldots, y_k, y_{k+1}, \ldots,
    y_{n}]^{S_k\times S_{n-k}}}
{(\prod (1+y_i))} \\
&=\dfrac{\Q[s_1, \ldots, s_k,q_1, \ldots, q_{n-k}]}
{((1+s_1+\cdots +s_k)(1+q_1+\cdots +q_{n-k}))},
\end{align*}
where $s_r= c_r(S)=\sigma_r(y_1, \ldots, y_k)$ and
$q_r= c_r(Q) = \sigma_r(y_{k+1}, \ldots, y_n)$.

\begin{prop} \label{p:charnogkn}
The characteristic numbers of 
$G(k,n)$ are
\begin{equation} \label{e:charnogkn}
\int_{G(k,n)} c_1(S)^{m_1} \cdots c_k(S)^{m_k} =
\sum_I \dfrac{\prod_{r=1}^k \gs_r (u_{i_1}, \dots, 
u_{i_k})^{m_r}}
{\prod_{i\in I} \prod_{j\in J} (u_i - u_j)},
\end{equation}
where $\sum m_r = k(n-k)$, $I$ runs over all multi-indices
$1 \le i_1 < \dots < i_k \le n$ and $J$ is its complementary multi-index.
\end{prop}

\begin{proof}
In Theorem \ref{t:charnogh}, take $f(y)$ to be $\prod_{r=1}^k c_r(S)^{m_r} = 
\prod_{r=1}^k
\gs_r(y_1, \dots, y_k)^{m_r}$
and 
$$w=(i_1, \ldots, i_k, j_1, \ldots, j_{n-k})$$
as in \eqref{e:permutation}.
Because $w\cdot u_1 = u_{i_1}$, \ldots, $w\cdot u_k = u_{i_k}$,
\[
w\cdot f(u) = \prod_{r=1}^k \gs_r(u_{i_1}, \dots, 
u_{i_k})^{m_r}.
\]
By \eqref{e:chernclass},
\[
\prod_{\ga \in \tri^+ \setminus \tri^+(H)} c_1(S_{\ga}) = 
\prod_{i=1}^k \prod_{j=k+1}^n (u_i - u_j).
\]
By \eqref{e:permutation},
\[
w\cdot \left(
\prod_{\ga \in \tri^+ \setminus \tri^+(H)} c_1(S_{\ga})
\right)=
w\cdot \left(
\prod_{i=1}^k \prod_{j=k+1}^n (u_i - u_j)\right) =
\prod_{i\in I} \prod_{j \in J} (u_i - u_j).    \tag*{\qed}
\]
\end{proof}

One of the surprising features of the localization formula 
is that although the right-hand side of \eqref{e:charnogkn} 
is apparently a sum of rational functions of $u_1, \dots, 
u_n$, the sum is in fact an integer.

\begin{exam*}
As an example, we compute the characteristic numbers of $\C
P^2=G(1,3)$.
The rational cohomology of $\C P^2$ is $H^*(\C P^2)= \Q[x]/(x^3)$,
generated by $x=c_1(S^{\spcheck}) = -c_1(S)$.
By Proposition~\ref{p:charnogkn},
\begin{align*}
\int_{\C P^2} x^2 &= \int_{G(1,3)} c_1(S)^2 = \sum_{i=1}^3
\frac{u_i^2}{\prod_{j\ne i}(u_i - u_j)}\\
&=\frac{u_1^2}{(u_1-u_2)(u_1-u_3)} +\frac{u_2^2}{(u_2-u_1)(u_2-u_3)}
+\frac{u_3^2}{(u_3-u_1)(u_3-u_2)},
\end{align*}
which simplifies to 1, as expected.
\end{exam*}

\end{document}